\documentclass[11pt]{amsart} 

\usepackage{amssymb,amsmath,amscd,amsxtra,color}
\usepackage[margin=4cm]{geometry}
\usepackage[all]{xy}

\newtheorem{theorem}{Theorem}[section]
\newtheorem{lemma}[theorem]{Lemma}

\newtheorem{proposition}[theorem]{Proposition}  
\newtheorem{corollary}[theorem]{Corollary}

\theoremstyle{definition}

\newcommand{\beq}{\begin{equation}}
\newcommand{\eeq}{\end{equation}}
\newcommand{\beqa}{\begin{eqnarray}}
\newcommand{\eeqa}{\end{eqnarray}}
\newcommand{\beaa}{\begin{eqnarray*}}
\newcommand{\ben}{\begin{eqnarray*}}
\newcommand{\eaa}{\end{eqnarray*}}
\newcommand{\een}{\end{eqnarray*}}

\newcommand{\CC}{\mathbb{C}}

\begin{document} 

\title[K-theoretic GW invariants in genus 0 and integrable hierarchies]
{K-theoretic Gromov--Witten invariants in genus 0 and integrable hierarchies}
\author{Todor Milanov and Valentin Tonita}

\address{Kavli IPMU (WPI), UTIAS, The University of Tokyo, Kashiwa, Chiba 277-8583, Japan}
\email{todor.milanov@ipmu.jp}
\address{Humboldt Univerit\"at, Rudower Chaussee 25 Raum 1425, Berlin
  10099, Germany}
\email{valentin.tonita@hu-berlin.de}

\begin{abstract}
We prove that the genus-0 invariants in K-theoretic Gromov--Witten
theory are governed by an integrable hierarchy of hydrodynamic type. 
If the K-theoretic quantum product is semisimple, then we also prove
the completeness of the hierarchy.
\end{abstract}
\maketitle

\setcounter{tocdepth}{2}
\tableofcontents

\section{Introduction}

Let $X$ be a smooth complex projective variety and $K(X)=K^0(X;\mathbb{C})$ be the Grothendieck group of topological
vector bundles on $X$. Let us denote by $X_{0,n,d}$ the moduli stack of genus-0 stable maps
of degree $d\in H_2(X;\mathbb{Z})$ with $n$ marked points. The
operation that assigns to a point in the moduli space the 
cotangent line at the $i$-th marked point is functorial and it gives
rise to a line bundle $L_i$, while evaluation at the $i$-th marked
point gives rise to a map of Deligne--Mumford stacks
$\operatorname{ev}_i:X_{0,n,d}\to X$ known as the evaluation map. 
If $E_1,\dots,E_n\in K(X)$, then following Givental and Y. P. Lee (see
\cite{Gi1, Lee1}) we introduce the K-theoretic Gromov--Witten
invariant
\ben
\langle E_1 L_1^{k_1},\dots,E_n L_n^{k_n}\rangle_{g,n,d} = 
\chi\Big(\mathcal{O}_{\rm virt}\otimes  
\operatorname{ev}^*_1(E_1) L_1^{k_1} \cdots \operatorname{ev}^*_n(E_n) L_n^{k_n} \Big),
\een
where $\chi(\mathcal{F})$ denotes the holomorphic Euler characteristic
of $\mathcal{F}$ and $\mathcal{O}_{\rm virt}$ is the so called {\em virtual structure
  sheaf} (see \cite{Lee1}). 

Let us fix a set $P_1,\dots,P_r$ of ample line bundles, s.t.,
$p_i=c_1(P_i)$ form a basis of $H_2(X;\mathbb{Z})/{\rm torsion}$. If
$d\in H_2(X;\mathbb{Z})$ then we define
\ben
Q^d:=Q_1^{\langle p_1,d\rangle}\cdots Q_r^{\langle p_r,d\rangle},
\een
where $Q_1,\dots,Q_r$ are formal variables known as the {\em Novikov
  variables}. For $t\in K(X)$ put
\beq\label{KGW-series}
\langle E_1 L_1^{k_1},\dots,E_n L_n^{k_n}\rangle_{g,n}(t) =
\sum_{d} \sum_{\ell=0}^\infty 
\frac{Q^d}{\ell !} \langle E_1 L_1^{k_1},\dots,E_n
L_n^{k_n},t,\dots,t\rangle_{g,n+\ell,d}.
\eeq
\subsection{Frobenius-like structure and the $S$-matrix}
Let us fix a basis $\{\Phi_i\}_{i=0}^N \subset K(X)$ and write
$t=\sum_{i=0}^N t_i\Phi_i$. In general, the RHS of
\eqref{KGW-series} is a formal power series in $Q_i$ $(1\leq i\leq r)$ and
$t_j$ $(0\leq j\leq N)$. 
We are going to assume that there exists an open subset 
\ben
B=\{(Q,t)\ |\ |Q_i|<\epsilon, |t_j|<\epsilon\}\subset
\mathbb{C}^r\times \CC^{N+1},
\een  
in which the series \eqref{KGW-series} is convergent for all possible
choices of the insertions $E_iL_i^{k_i}$. This condition is
satisfied in many important cases, such as 
when $X$ is a projective space or a flag manifold (see \cite{IMT}). 
Our results are true in general as
well, but we have to replace $B$ with the formal germ of
$\mathbb{C}^r\times \CC^{N+1}$  at $(0,0).$
Slightly abusing the notation we denote by $TB$ the relative tangent
bundle of the projection $B\to \CC^r$, $(Q,t)\mapsto Q$. 

The vector space $K(X)$ is equipped with the Euler pairing
\ben
g(\Phi_i\otimes \Phi_j) :=  
g_{ij}:=
\int_X \operatorname{ch}(\Phi_i) \operatorname{ch}(\Phi_j) \operatorname{Td}(X).  
\een
Put 
\ben
G(\Phi_i\otimes \Phi_j) := 
G_{ij}(t) :=
g_{ij}+\langle \Phi_i,\Phi_j\rangle_{0,2}(t).
\een
The $K$-theoretic quantum product is defined by 
\ben
G(\Phi_i\bullet_t \Phi_j, \Phi_k) = \langle\Phi_i,\Phi_j,\Phi_k\rangle_{0,3}(t).
\een
Using the natural trivialization $TB=B\times K(X)$, $\partial/\partial
t_i\mapsto \Phi_i$, we can identify $g$ and $G$ with non-degenerate
bilinear forms on $TB$ and $\bullet$ with a multiplication on $TB$.
It is known that the multiplication $\bullet$ and the pairing $G$
satisfy Dubrovin's axioms of a Frobenius manifold (see \cite{Du})
except for flat identity
and Euler vector field (see \cite{Gi1}).

Following Givental \cite{Gi1}, we introduce the operator series $S\in
\operatorname{End}(K(X))\otimes \mathbb{C}(q)[\![Q,t]\!]$ defined by  
\ben
G(\Phi_i,S(t,q)\Phi_j) = S_{ij}(t,q^{-1}),
\een
where 
\ben
S_{ij}(t,q):= g_{ij}+\left\langle \Phi_i,\frac{\Phi_j}{ 1-qL}\right\rangle_{0,2}(t).
\een 
Here $q$ is a complex number such that $(1-qL)$ is invertible in
$K(X_{0,n,d})$ for all $n$ and $d$. It is known that the
coefficient in front of a fixed monomial in $Q$ and $t$ of
$S_{ij}(t,q)$ is a rational function in $q$ which vanishes at
$q=\infty$ and has poles only at the roots of unity. The operator series
$S(t,q)$ is a fundamental solution to the following system of
differential equations
\ben
(q-1)\partial_{t_i} S(t,q) = \Phi_i\bullet_t S(t,q),\quad 0\leq i\leq N.
\een
In particular, the formulas
\ben
\nabla_{\partial/\partial t_i} := \partial_{t_i} -(q-1)^{-1}
\Phi_i\bullet_t,\quad 0\leq i\leq N,
\een
define a family of flat connections on $TB$. It is easy to check that
for $q=-1$ the connection $\nabla$ turns into the
Levi--Cevita connection for the pairing $G$, which in particular
proves that the bilinear pairing $G$ is  flat.

\medskip

Motivated by  Dubrovin's construction of the principal hierarchies in the theory of
Frobenius manifolds (see\cite{Du, DuZ}), we introduce a
hierarchy of differential equations on the formal loop space
$\mathcal{L}(B)$. It is convenient however to work first 
purely algebraically using the language of differential algebras and
then reformulate the construction using the geometry of the formal
loop space.  

\subsection{Commuting flows}

Let $\CC[Q,v,\partial v,\partial^2v,\dots]$ be the free differential
algebra over the ring $\CC[Q]$ on the vector variable
$v=(v_0,\dots,v_N)$.  We refer to
$\partial^kv=(\partial^kv_0,\dots,\partial^kv_N)$ with $k>0$ as {\em jet
  variables}. Each $\partial^kv_i$ should be viewed as a formal
variable, $\partial^0v_i:=v_i$, and the differentiation is defined in the obvious way
$\partial(\partial^kv_i)=\partial^{k+1}v_i$. Let us denote by
$\mathcal{A}$ the completion of the differential algebra 
consisting of polynomials in the jet variables whose coefficients are
convergent power series for $(Q,v)\in B$.  If $q\in
\CC\setminus{\{0\}}$ is not a root of unity, then we define 
derivations $T_i(q)$ $(0\leq i\leq N)$ of $\mathcal{A}$, such that $T_i(q)$ commutes with $\partial$
and 
\ben
T_i(q) v := \partial(S(v,q)\Phi_i  )=(q-1)^{-1}\partial v\bullet S(v,q)\Phi_i,
\een 
where both sides take value in
$\mathcal{A}\otimes K(X)$ via the identification $v:=\sum_{i=0}^N v_i
\Phi_i$, and the derivation $\partial$ of $\mathcal{A}$ is extended
uniquely to a derivation on $\mathcal{A}\otimes K(X)$, such that
$\partial(1\otimes \Phi_i)=0$.

Following Dubrovin and Novikov (see \cite{DuN}) we equip the vector space
$\mathcal{F}:=\mathcal{A}/\partial \mathcal{A}$ with a Poisson bracket
of {\em hydrodynamic type}
\ben
\{f_1,f_2\} =\sum_{i,j=0}^N\int 
\frac{\delta f_1}{\delta v_i} A^{ij} \frac{\delta f_2}{\delta v_j} ,
\quad f_1,f_2\in \mathcal{F},
\een
where $\int:\mathcal{A}\to \mathcal{F}$ is the projection
map, $A^{ij}$ is the first order differential operator
\ben
A^{ij} = G^{ij}(v) \partial -\sum_{s,k=0}^N
G^{is}\Gamma_{sk}^j(v) \partial v_k,
\een
$G^{ij}$ are the entries of the inverse of the matrix $(G_{ij})$, 
\ben
\Gamma_{sk}^j(v)=
\frac{1}{2}\sum_{\ell=0}^N G^{j\ell}\, G(\Phi_s\bullet\Phi_k,\Phi_\ell)
\een 
are the Christoffel's symbols of the pairing $G$,
and 
\ben
\frac{\delta}{\delta v_i}=\sum_{\ell=0}^\infty (-\partial)^\ell
\frac{\partial}{\partial(\partial^\ell v_i)}
\een
is the variational derivative. Recall that $\delta/\delta
v_i\circ \partial =0$ (in $\mathcal{A}$), so the variational
derivative descends to the quotient $\mathcal{F}$.

Given $H\in \mathcal{F}$, we can define a derivation $X_H$ of
$\mathcal{A}$, such that it commutes with $\partial$ and
\ben
X_H(v_i)=\sum_{j=0}^N \Big(G^{ij}(v)\partial - \sum_{s,k=0}^N
G^{is}(v)\Gamma_{sk}^j(v) \partial v_k\Big) \, \frac{\delta H}{\delta v_j}.
\een
The map $H\to X_H$ turns $\mathcal{A}$ into a Poisson
$\mathcal{F}$-module:
\ben
[X_{H_1},X_{H_2}] = -X_{\{H_1,H_2\} },\quad H_1,H_2\in \mathcal{F}.
\een
We will refer to $H$ as {\em Hamiltonian} and to
$X_H$ as {\em Hamiltonian derivation}. We say that the Hamiltonian
derivation $X_H$ is polynomial in $v_0$ if the differential polynomial
$X_H(v)$ depends polynomially in $v_0$. Our first result can be stated as
follows.
\begin{theorem}\label{t1}
a) The derivations $T_i(q)$ ($i=0,1,\dots,N$) are Hamiltonian and
$T_i(q_1)$ and $T_j(q_2)$ commute for all $i$, $j$, $q_1$, and $q_2$. 

b) If the K-theoretic quantum product is semisimple, then every Hamiltonian
derivation $X_H$ commuting with $T_i(q)$ for all $i$ and polynomial in
$v_0$ has the form 
\ben
X_H(v) = -\operatorname{Res}_{q=\infty} \partial S(v,q) C(q)dq,
\quad C(q)\in K(X)[q].
\een
\end{theorem}
Let us point out that the Hamiltonian corresponding to the derivation
$T_i(q)$ is 
\beq\label{H_i}
H_i(v,q) := \frac{2}{1+q}\Big(g(v,\Phi_i)+\Big\langle
  \frac{\Phi_i}{1-q^{-1}L}\Big\rangle_{0,1}(v) \Big)
\quad (\mbox{mod }\partial\mathcal{A}).
\eeq
Part a) follows immediately from the fact that $S(v,q)$ is a
symplectic transformation solving the K-theoretic quantum differential
equations, while for part b) we follow the argument of the proof of
Lemma 3.6.7 and Lemma 3.6.20 in \cite{DuZ}.  
 
\subsection{Integrable hierarchies}
Theorem \ref{t1} allows us to construct an integrable hierarchy of
evolutionary PDE on the loop space $\mathcal{L}(B):=C^\infty(S^1,B).$ 
Let us parametrize $S^1$ via $x\mapsto e^{2\pi\sqrt{-1} x}$.
Recall that a {\em local functional} is a function on $\mathcal{L}(B)$
of the form
\ben
v(x)\mapsto \int_0^1  P(x,v(x),\partial_xv(x),\dots)dx,
\een
where 
\ben
P(x,v,\partial v,\partial^2 v,\dots)=
\sum_{\substack{k=(k_1,\dots,k_s):k_a\geq 1\\i=(i_1,\dots,i_s)}} 
c_{k,i}(x,v) 
(\partial^{k_1}v_{i_1}) \cdots (\partial^{k_s}v_{i_s}), 
\een
is a polynomial in the jet variables whose coefficients
$c_{k,i}(x,v)\in C^\infty(S^1\times B)$. In particular, the elements
of $\mathcal{F}$ are local functionals and the
definition of the Poisson bracket of hydrodynamic type extends to the
entire space of local functionals. The loop space $\mathcal{L}(B)$ has
also another class of functions that comes from the evaluation map
\ben
\operatorname{ev}:S^1\times \mathcal{A}\to C^\infty(\mathcal{L}(B))
\een
defined by 
\ben
\operatorname{ev}(x,P)(v)=P(v(x),\partial_x v(x),\partial_x^2v(x),\dots).
\een 
Following the notation in physics we write the evaluation maps as
\ben
\operatorname{ev}(x,P)(v) =
\int_0^1 
P(v(y),\partial_y v(y),\partial_y^2v(y),\dots)
\delta(x-y) dy,
\een
which allows us to think of $\operatorname{ev}(x,P) $ as a local
functional and to extend the definition of the Poisson bracket. Note
that we have 
\ben
\{ \operatorname{ev}(x,P),H\} = \operatorname{ev}(x, X_H(P)),\quad
P\in \mathcal{A},\quad H\in \mathcal{F},
\een
which justifies why the derivation $X_H$ was called Hamiltonian. 
Put
\ben
P_{n,i}(v,\partial v) := -\operatorname{Res}_{q=\infty} dq\, (q-1)^{n}
T_i(q)v,\quad 0\leq i\leq N, \quad n\geq 0.
\een
Let us introduce the following system of ODEs  
\beq\label{is-alpha}
\partial_{t_{n,i}} v= P_{n,i}(v,\partial_x v), 
\quad 0\leq i\leq N, \quad n\geq 0,
\eeq
where
\ben
v=v(x,\mathbf{t}),\quad \mathbf{t}=(t_{n,i})^{i=0,\dots,N}_{n=0,1,2,\cdots}.
\een
Note that the first equation of the hierarchy is $\partial_{0,0}v
=\partial_x v$. Therefore, we may assume that $x=t_{0,0}$.
\begin{corollary}\label{c1}
The equations \eqref{is-alpha} are Hamiltonian and pairwise
compatible.  If the quantum K-product is semisimple, then
\eqref{is-alpha} is a completely integrable system. 
\end{corollary}
The Hamiltonians corresponding to \eqref{is-alpha} are given by 
\ben
H_{n,i}(v):=-\operatorname{Res}_{q=\infty} dq\, (q-1)^{n} H_i(v,q).
\een
The Taylor's series expansion of $S(v,q)$ near $q=\infty$ has the form 
\ben
S(v,q)=1+\sum_{n=1}^\infty S_n(v) (q-1)^{-n}.
\een
The quantum differential equations imply that $\partial_{v_0}S_n(v) =
S_{n+1}(v)$. From this relation we get that $S_n(v)$ depends
polynomially on $v_0$. By definition
\ben
X_{H_{n,i}}(v) = P_{n,i}(v,\partial v)= -\operatorname{Res}_{q=\infty} dq\, (q-1)^{n-1}
\partial v\bullet S(v,q)\Phi_i = \partial v\bullet S_n(v)\Phi_i,
\een
so the Hamiltonian derivations $X_{H_{n,i}}(v)$ are polynomial in
$v_0$. Corollary \ref{c1} is a direct consequence of 
Theorem \ref{t1}. 

\subsection{The topological solution}

The second goal of our paper is to construct a solution to the 
system \eqref{is-alpha} in terms of genus-0 K-theoretic GW invariants. The genus-0
total descendant potential is defined by 
\ben
\mathcal{F}(\mathbf{t}) = 
\sum_{n=0}^\infty 
\frac{1}{n!}
\langle \mathbf{t}(L_1),\dots,\mathbf{t}(L_n)\rangle_{0,n},
\een
where $\mathbf{t}(q)\in K(X)[q,q^{-1}]$. Let us define 
\ben
\mathbf{t}(q):=\sum_{k=0}^\infty \sum_{i=0}^N t_{k,i} \Phi_i (q-1)^k,
\een
where $t_{k,i}$ are formal variables.
Let us denote $\partial_{m,i}:=\partial/\partial t_{m,i}$. 
Recall the K-theoretic J-function of
$X$
\ben
J(v,q) = 1-q+v+\sum_{i=0}^N\Phi^i\langle \frac{\Phi_i}{1-qL}\rangle_{0,1}(v)=(1-q)S(v,q)^{-1}\mathbf{1},
\een
where $v\in K(X)$ and $\Phi^i:=\sum_{j=0}^N g^{ij}\Phi_j$ ($0\leq i\leq N$)  is the basis of $K(X)$ dual
to $\{\Phi_i\}_{i=0}^N$  with respect to the Euler pairing $g$.  Using
the $J$-function we introduce new local coordinates 
$w=\sum_{a=0}^N w_a\Phi^a$ on $B$
\beq\label{Miura-tr}
w=J(v,0)-1 = v+\sum_{i=0}^N \langle\Phi_i\rangle_{0,1}(v) \Phi^i.
\eeq
Note that the vector fields $\partial/\partial w_a =\sum_b
G^{ab}(v)\partial/\partial v_b$ correspond to the differential forms
$dv_a$ under the identification $TB\cong T^*B$ defined via the
 bilinear pairing $G$. 
\begin{theorem}\label{t2} 
The functions
\ben
w_a(\mathbf{t})= \partial_{0,a}\partial_{0,0}\mathcal{F}(\mathbf{t}),\quad 0\leq a\leq N,
\een
provide a solution to the integrable hierarchy obtained from
\eqref{is-alpha}  via the Miura transformation \eqref{Miura-tr}. 
\end{theorem}

\subsection{Example: $X=pt$}
In this case, the quantum differential equations are straightforward
to solve. We have
\ben
S(v,q)=e^{v/(q-1)},\quad
J(v,q)=e^{v/(1-q)}.
\een
The integrable system \eqref{is-alpha} takes the form
\ben
\partial_{n} v = \frac{v^n}{n!}\,\partial v,\quad n=0,1,\dots,
\een
i.e., this is the dispersionless KdV hierarchy. 
The Miura transformation is $w=e^v-1$, so the dispersionless KdV
hierarchy is transformed into
\ben
\partial_n w = \frac{1}{n!}(\log(1+w))^n\, \partial w,\quad n=0,1,\dots.
\een
The topological solution is
\ben
w(\mathbf{t}) = 
\sum_{n=1}^\infty \frac{1}{n!} 
\langle 1,1,\mathbf{t}(L),\dots,\mathbf{t}(L)\rangle_{0,n+2},
\een
where
\ben
\mathbf{t}(L) = \sum_{n=0}^\infty t_n (L-1)^n.
\een
Using the above equation and the string equation we can compute all
genus 0 K-theoretic GW invariants of the point. Our answer agrees with
the answer of Y.P. Lee \cite{Lee2}.
It is very tempting to conjecture that the higher genus deformation of
the topological solution
\ben
\widehat{w}(\hbar,\mathbf{t}):= 
\sum_{g,n=0}^\infty \frac{\hbar^g}{n!} 
\langle 1,1,\mathbf{t}(L),\dots,\mathbf{t}(L)\rangle_{g,n+2},
\een
is a solution to an integrable hierarchy which is a Miura transform of
the KdV hierarchy. 

While the higher-genus K-theoretic GW invariants of the point are still very
difficult to compute, motivated by the above example, we would like to investigate
further the relation of the hierarchy \eqref{is-alpha} with the
genus-0 hierarchies in cohomological GW theory. Especially, the adelic
characterization of genus 0 K-theoretic GW invariants  (see \cite{GT})
seems to be appropriate to study such a question. In particular, it
would be interesting to find out if the integrable hierarchy
\eqref{is-alpha} is bi-Hamiltonian. We are planning to investigate
these questions in a future publication. 

\subsection{Acknowledgements}
We are thankful to Y. Zhang for showing interest in our work and for
several useful discussions. The work of T.~M.~is partially supported by JSPS Grant-In-Aid 26800003 
and by the World Premier International Research Center Initiative (WPI
Initiative),  MEXT, Japan.  This work started during a visit 
of the second author to Kavli-IPMU. V.~T. would like to thank  the institute for hospitality and support.

\section{Proof of Theorem \ref{t1}}
The proof of part a) will be given in Sections \ref{ham} and
\ref{comm}. The argument for part b) will be given in Section
\ref{compl}.

\subsection{The derivations are Hamiltonian}\label{ham}
The fact that the derivation $T_i(q)$ is Hamiltonian with
Hamiltonian \eqref{H_i} follows from the following formula:
\beq\label{bracket:H_i}
\{f,H_i\}=\int \sum_{a,s=0}^N\,
G(\partial S(v,q)\Phi_i,\Phi_s)\, G^{sa}\,
\frac{\delta f}{\delta v_a}.
\eeq
Therefore, we just need to prove \eqref{bracket:H_i}. 

Let us denote by $g=(g_{ab})$ and
$G=(G_{ab})$ the matrices of the two pairings. The pairing $G$ can be expressed in terms of $g$ and
$S(v,0)$. Indeed, we have the following identity expressing the
fact that $S$ is a symplectomorphism (see \cite{Gi1, IMT})
\beq\label{S:sympl}
g(a,S(v,q^{-1})^{-1}b) = G(S(v,q)a,b), \quad a,b\in K(X).
\eeq
If we set $q=0$, then since $S(v,\infty)=1$, we get
\ben
g=G\, S(v,0).
\een
On the other hand,
\ben
\frac{1}{2}(q+1)A^{ab} \frac{\delta H_i}{\delta v_b} = G^{ab}\partial \,
G(S(v,q)\Phi_i,\Phi_b) -\sum_{s,\ell} G^{as}\Gamma_{s\ell}^b(v) \partial v_\ell G(S(v,q)\Phi_i,\Phi_b).
\een
Note that $\partial G^{-1}=-\sum_\ell \Omega_\ell G^{-1}\partial
v_\ell$, where $\Omega_\ell=\Phi_\ell\bullet$ is the operator of
quantum K-theoretic multiplication by $\Phi_\ell$. Therefore,
\ben
G^{ab}\circ\partial =\partial\circ G^{ab} +\sum_{s,\ell}
(\Omega_\ell)_{as}G^{sb}\partial v_\ell.
\een
By definition 
\ben
\sum_{k,\ell}
(\Omega_\ell)_{ak}G^{kb}\partial v_\ell=2 \sum_{s,\ell} G^{as}\Gamma_{s\ell}^b(v) \partial v_\ell,
\een
so
\ben
\frac{1}{2}(q+1)A^{ab} \frac{\delta H_i}{\delta v_b} =
\partial \,\Big( G^{ab}
G(S(v,q)\Phi_i,\Phi_b)\Big)+
\frac{1}{2} \sum_{s,\ell}
(\Omega_\ell)_{as}G^{sb}\partial v_\ell G(S(v,q)\Phi_i,\Phi_b).
\een
Note that 
\ben
\sum_b (\Omega_\ell)_{as}G^{sb}G(S(v,q)\Phi_i,\Phi_b) = 
G(S(v,q)\Phi_i,\Phi_\ell\bullet \widetilde{\Phi}^a)=
G(\Phi_\ell\bullet  S(v,q)\Phi_i, \widetilde{\Phi}^a),
\een
where $\widetilde{\Phi}^a=\sum_{b=0}^N G^{ab}\Phi_b$.
Recalling the differential equation
$(q-1)\partial_{v_\ell}S(v,q)=\Phi_\ell\bullet S(v,q)$ we get
\ben
\frac{1}{2}(q+1)\sum_b A^{ab} \frac{\delta H_i}{\delta v_b} =
\frac{1}{2}(q+1)\sum_b G^{ab}
G(\partial S(v,q)\Phi_i,\Phi_b),
\een
where we used that 
\ben
\sum_b \partial \Big(G^{ab}G(x ,\Phi_b)\Big) = \sum_b \partial \Big(g^{ab}
g(x,\Phi_b) \Big)= \sum_b G^{ab}G(\partial x,\Phi_b).
\een
By definition 
$\{f,H_i\} = \int \sum_{a,b} \frac{\delta f}{\delta v_a}
A^{ab} \frac{\delta H_i}{\delta v_b}
$, so the above formula implies \eqref{bracket:H_i}.

\subsection{The Hamiltonians are in involution}\label{comm}

We have to prove that 
\ben
\{H_i(v,q_1),H_j(v,q_2)\}=0,\quad  \forall
i,j=0,\dots,N.
\een 
Recalling formula \eqref{bracket:H_i} we get that the vanishing of the
above bracket is equivalent to the fact that 
the differential polynomial
\ben
G(\partial S(v,q)\Phi_i,S(v,q_2)\Phi_j) = 
\frac{1}{q_1-1}
\sum_{\ell=0}^N g(S(v,q_2^{-1})^{-1} \Omega_\ell
S(v,q_1)\Phi_i,\Phi_j)\partial v_\ell
\een
is $\partial$-exact, where $\Omega_\ell$ is the linear operator of
quantum multiplication by $\Phi_\ell$ and the equality is derived by
using the quantum differential equations for $S(v,q_1)$ and the
symplectic property of $S(v,q_2)$ (see \eqref{S:sympl}). Recalling
again the quantum differential equations for $S$ we get 
\ben
S(v,q_2^{-1})^{-1} \Omega_\ell S(v,q_1) = 
\partial_{v_\ell} \left(
\frac{(q_1-1)(q_2-1)}{q_1q_2-1} 
\Big(S(v,q_2^{-1})^{-1}S(v,q_1)-1\Big)\right).
\een
The above formula implies the $\partial$-exactness, because the
pairing $g$ is constant.

\subsection{Canonical coordinates} The notion of canonical coordinates
of a semisimple Frobenius manifold (see \cite{Du}) is straightforward to generalize to
the settings of quantum K-theory. In fact, the quantum K-ring is a
F-manifold in the sense of Hertling and Manin (see \cite{HM}), so we
can recall the notion of canonical coordinates of a  F-manifold (see
Corollary 5.2.1 in \cite{Ma}). However, for the sake of completeness, let us
work out the necessary details in our setting.

We say that the K-theoretic quantum product is {\em semisimple} at
$v^{(0)}\in B$ if there are local holomorphic vector fields $e_i$, $0\leq i\leq N$,
defined in an open neighborhood of $v^{(0)}$, such that 
\ben
e_i\bullet_v  e_j=\delta_{ij} e_j,\quad i,j=0,\dots N,
\een
for all $v$ sufficiently close to $v^{(0)}$. Note that in the basis
$\{e_i\}_{i=0}^N$ the pairing $G$ is diagonal 
\ben 
G(e_i,e_j)=\Delta_i(v) \delta_{ij}
\een
where $\Delta_i(v)\neq 0$ for all $i$ and for all $v$ sufficiently
close to $v^{(0)}$. 
\begin{lemma}\label{can_coord}
If $v^{(0)}\in B$ is a semisimple point and $\{e_i\}_{i=0}^N$ is a corresponding
basis of idempotents, then there are local coordinates $u_i$, $0\leq
i\leq N$, such that $e_i=\partial/\partial u_i$. 
\end{lemma}
\proof
We need to prove that the vector fields $e_i$ and $e_j$ commute. 
The quantum connection can be written in the form
\ben
\nabla = \nabla^{\rm L.C.} + \alpha \sum_{\ell=0}^N \Omega_\ell(v)
dv_\ell,
\quad
\alpha := -\frac{1}{2}-\frac{1}{q-1},
\een
where $\nabla^{\rm L.C.}$ is the Levi--Civita connection for the
bilinear pairing $G$. The flatness of $\nabla$ implies that 
\ben
(\nabla_{e_i}\nabla_{e_j}-\nabla_{e_j}\nabla_{e_i} )e_k=
\nabla_{[e_i,e_j]}e_k. 
\een
Comparing the coefficients in front of $\alpha$ we get 
\ben
e_i\bullet (\nabla^{\rm L.C.}_{e_j} e_k)+\nabla^{\rm L.C.}
_{e_i}(e_j\bullet e_k) - (i\leftrightarrow j) = [e_i,e_j]\bullet e_k,
\een
where the expression in the brackets on the LHS is obtained from the
preceeding expression by switching $i$ and $j$. Using that $e_a$ are
idempotents we get 
\ben
\Gamma_{jk}^i e_i - \Gamma_{ik}^j e_j +\sum_{s=0}^N
(\delta_{jk}-\delta_{ik})\Gamma_{ij}^s e_s = [e_i,e_j]\bullet e_k.
\een
Multiplying both sides by $e_k$, we get that the LHS is $0$, while the
RHS does not change, so $[e_i,e_j]\bullet e_k=0$ for all $k$ and we
get $[e_i,e_j]=0$.
\qed

Let us denote by $\Psi$ the matrix with entries $\Psi_{ai}
=\frac{\partial v_a}{\partial u_i}$, then the above lemma implies that 
\ben
\Psi^{-1}\Omega_\ell \Psi = \operatorname{Diag}
\Big( 
\frac{\partial u_0}{\partial v_\ell},\dots,\frac{\partial u_N}{\partial v_\ell}
\Big).
\een

\subsection{Completeness}\label{compl}
Let us assume that $f\in \mathcal{A}$ is such that $\{\int f,H_i(v,q)\}=0$
for all $i=0,\dots,N$. 
\begin{lemma}\label{jet_indep}
The differential polynomial $f$ is independent of the jet variables,
i.e., $\partial f/\partial (\partial^kv_a)=0$ for $k>0$.
\end{lemma}
\proof
Recalling formula \eqref{bracket:H_i} and using also that 
\ben
\sum_{b=0}^N G^{ab}G(x,\Phi_b) = \sum_{b=0}^N g^{ab}g(x,\Phi_b)
\een
we get that the differential polynomial 
\ben
\sum_{j=0}^N W^j(v,\partial v,\dots)\partial g(S(v,q)\Phi_i,\Phi_j)
\een
is $\partial$-exact, where $W^j=\sum_a g^{ja}\frac{\delta f}{\delta
  v_a}$. Therefore the variational derivatives 
\ben
\frac{\delta}{\delta v_a} 
\Big(\sum_{j=0}^N W^j(v,\partial v,\dots)\partial
g(S(v,q)\Phi_i,\Phi_j)\Big) =0
\een
for all $a=0,\dots,N$. By definition, the above variational derivative
is the sum of
\beq\label{var_der:1}
\partial_{v_a} \Big(\sum_j W^j \partial
g(S(v,q)\Phi_i,\Phi_j)\Big) -\partial \partial_{\partial v_a} \Big(\sum_j W^j \partial
g(S(v,q)\Phi_i,\Phi_j)\Big)
\eeq
and 
\beq\label{var_der:2}
\sum_{k\geq 2}\sum_{j=0}^N
(-\partial)^k \Big(\frac{\partial W^j}{\partial (\partial^k v_a)} \partial
g(S(v,q)\Phi_i,\Phi_j)\Big).
\eeq
Let us denote by $\Omega_\ell$ the linear operator of quantum
multipication by $\Phi_\ell$. Recalling also the quantum differential
equation for $S$ and canceling a factor of $q-1$, we transform \eqref{var_der:1} and
\eqref{var_der:2} respectively into 
\ben
g\Big(
\sum_{j,\ell} \Big(\partial_{v_a}(W^j\Omega_\ell) S-
\partial\Big(
\partial_{\partial v_a} W^j\Omega_\ell S\Big)\Big)
\Phi_i,\Phi_j\Big) 
\partial v_\ell
- \sum_{j=0}^N \partial(W^j\Omega_a)S\Phi_i,\Phi_j\Big)
\een
and 
\ben
\sum_{k\geq 2}
\sum_{k'=0}^k {k\choose k'}
\sum_{\ell=0}^N
g\Big(
(-\partial)^{k-k'}\Big(
\frac{\partial W^j}{\partial (\partial^k v_a)} \Omega_\ell \partial v_\ell\Big)
(-\partial)^{k'}S\Phi_i,\Phi_j)\Big).
\een
The sum of the above expressions can be written in the form
$\sum_jg(A_jS\Phi_i,\Phi_j)=0$, where $A_j$ is a polynomial in
$(q-1)^{-1}$ with coefficients in $\operatorname{End}(K(X))$. Since
$g$ is a non-degenerate pairing, we can cancel $S$, i.e. we must have
$\sum_{j=0}^Ng(A_j\Phi_i,\Phi_j)=0$. Let $m$ be a maximal non-negative integer,
s.t., $\partial_{\partial^m v_a}W^j\neq 0$ for some $j$. If $m\geq 2$,
then the coefficient in front of $(q-1)^{-m}$ of $A_j$ is precisely 
\ben
(-1)^m \frac{\partial W^j}{\partial (\partial^m v_a)} 
\Big(\sum_{\ell=0}^N \Omega_\ell \partial v_\ell\Big)^{m+1}.
\een
If the quantum multiplication is semisimple, then we can find a
matrix $\Psi=\Psi(v)$, holomorphically invertible for all $v$ in a
neighborhood of  a fixed semisimple point, such that
\ben
\Psi^{-1}\Omega_\ell \Psi = \operatorname{Diag}\Big(\frac{\partial
  u_0}{\partial v_\ell},\dots, \frac{\partial u_N}{\partial v_\ell} \Big),
\een
where $u_i=u_i(v)$ $(0\leq i\leq N)$ are the corresponding canonical
coordinates. Using again that $g$ is non-degenerate,
we may replace $\Phi_i$ by $\Psi e_i$, where $e_i$ is the $i$th
standard vector. We get 
\beq\label{contrad}
-(-\partial u_i)^{m+1} \sum_{j=0}^N \frac{\partial W^j}{\partial
  (\partial^m v_a)}  g(\Psi e_i,\Phi_j) = 0.
\eeq
Let us choose a non-constant smooth loop $v:S^1\to B$, s.t., $v(x)$ is in a
neighborhood of the fixed semisimple point for all $x\in S^1$ and 
\ben
\partial_x u_i(v(x))\neq 0\ (0\leq i\leq N), \quad 
\operatorname{ev}\Big(x,\frac{\partial W^j}{\partial
  (\partial^m v_a)}\Big)(v)\neq 0
\een
for some $x\in S^1$ and some $j$. Such a choice is possible, because
otherwise either $u_i$ $(0\leq i\leq N)$ will fail to be local coordinates or we will
have a contradiction with the choice of $m$. However, evaluating
\eqref{contrad} at $v(x)$ yields a contradiction with the non-degeneracy
of $g$. Therefore, $m\leq 1$. A similar 
argument can be used to prove that $m$ can not be 1.
\qed

\begin{lemma}\label{comm_Ham:eqn}
If $f(v)\in \mathcal{A}$ is independent of the jet variables, then
vanishing of the brackets
$$
\Big\{\int f(v),H_i(v,q)\Big\}=0,
\quad \forall i=0,\dots,N
$$ 
is equivalent to the following system of equations 
\ben
\frac{\partial^2 f}{\partial v_i \partial v_a} = 
\sum_{k=0}^N
(\Omega_i)_{ka} 
\frac{\partial^2 f}{\partial v_k \partial v_0},\quad \forall i,a=0,\dots,N.
\een
\end{lemma}
\proof
Let us use the same notation as in the proof of Lemma
\ref{jet_indep}. Using that $W^j$ is independent of the jet variables,
then the vanishing of the variational derivatives, i.e. the sum of
\eqref{var_der:1} and \eqref{var_der:2} simplifies to
\ben
g\Big( 
\Big(
(\partial_{v_a} W^j)\Omega_\ell - 
(\partial_{v_\ell} W^j)\Omega_a
\Big)\Phi_i,
\Phi_j\Big)=0,
\een
where we used that $\partial_{v_a}\Omega_\ell
= \partial_{v_\ell}\Omega_a$, which follows from the flatness of the
quantum connection. Put $\ell=0$, so that $\Omega_\ell=1$. Then the
above identity becomes equivalent to the identity that we wanted to
prove. 
\qed

\emph{Proof of part b) of Theorem \ref{t1}}. According to
Lemma \ref{jet_indep} the Hamiltonian $H=f(v)$ is independent of the
jet variables. Let us denote by $F(v)$ the gradient of $f$ with
respect to the pairing $G$, i.e., 
\ben
G(F(v),\Phi_a) = \frac{\partial f}{\partial v_a}. 
\een
A straightforward computation, similar to the proof of formula
\eqref{bracket:H_i}, yields 
\beq\label{grad_de}
\frac{\partial F}{\partial v_\ell} = \Omega_\ell \, \frac{\partial
  F}{\partial v_0}, 
\quad 
\forall \ell=0,\dots,N,
\eeq
and
\beq\label{der:f}
X_H(v) = \partial F+\frac{1}{2}\partial v \bullet F.
\eeq
Note that the derivation $X_H$ is polynomial in $v_0$ if and only if
the gradient $F(v)$ depends polynomially in $v_0$. Arguing by
induction on the degree of $F$ as a polynomial in $v_0$ we prove that 
\beq\label{F:Res_formula}
F(v)=-\operatorname{Res}_{q=\infty} S(v,q)C(q)\frac{dq}{q-1}
\eeq
for some $C(q)\in K(X)[q]$. The above identity, formula
\eqref{der:f}, and the quantum differential equations 
$(q-1)\partial S(v,q) =\partial v\bullet S(v,q)$ imply the formula
stated in the Theorem. 

If the degree of $F$ is 0, then recalling the differential equation
\eqref{grad_de} we get  $\partial_{v_\ell}F=0$ for all $\ell$, so $F$ is
independent of $v$. We just have to use that 
\ben
1=-\operatorname{Res}_{q=\infty} S(v,q)\frac{dq}{q-1}.
\een
For the induction step, note that $\partial_{v_0} F$ is the gradient
of $-f+\partial_{v_0}f$. According to Lemma \ref{comm_Ham:eqn}, the
latter Poisson commutes with all $H_i(v,q)$, so using the inductive
assumption we have
\ben
\partial_{v_0} F(v)=-\operatorname{Res}_{q=\infty} S(v,q)C(q)\frac{dq}{q-1}
\een
for some $C(q)\in K(X)[q]$. Put 
\ben
\widetilde{F}(v) := -\operatorname{Res}_{q=\infty} \partial S(v,q)C(q)(q-1)\frac{dq}{q-1}.
\een
Recalling the differential equations \eqref{grad_de} we get
$\partial_{v_\ell}F=\partial_{v_\ell}\widetilde{F}$ for all $\ell$, hence
 $F-\widetilde{F}$ is independent of $v$. Since both
$\widetilde{F}$ and $F-\widetilde{F}$  can be presented in the form
\eqref{F:Res_formula}, the gradient $F$ also has such a presentation.
\qed

\section{Construction of the topological solution}
\subsection{From descendants to ancestors}
The genus-0 total ancestor potential $\overline{\mathcal{F}}_\tau$,
$\tau\in K(X)$, is defined by 
\ben
\overline{\mathcal{F}}_\tau(\mathbf{t})=\sum_{\ell,n=0}^\infty\sum_d
\frac{Q^d}{n!\ell!}\langle \mathbf{t}(\overline{L}_1),\dots,  \mathbf{t}(\overline{L}_n),\tau,\dots,\tau)\rangle_{0,n+\ell,d},
\een
where $\overline{L}_i$ ($1\leq i\leq n$) is the pullback of the tautological cotangent
line bundle $L_i$ on $\overline{\mathcal{M}}_{0,n}$via the forgetfull morphism
$X_{0,n+\ell,d}\to \overline{\mathcal{M}}_{0,n}$. By definition the
potential $\overline{\mathcal{F}}_\tau(\mathbf{t})$ is a formal power
series in $\mathbf{t}$, $\tau$, and $Q$. 
 Let us recall the following formula relating the ancestor and descendant potentials (see
 \cite{TT}, Appendix B)
\beq\label{F-anc}
\mathcal{F}(\mathbf{t}) = \frac{1}{2}\langle\mathbf{q}(L),\mathbf{q}(L)
\rangle_{0,2}(\tau) +\overline{\mathcal{F}}_\tau([S(\tau,q)\mathbf{q}(q)]_+).
\eeq
where $\mathbf{q}(q) := \mathbf{t}(q)+1-q$ and $[\ ]_+$ denotes the
operation induced from the projection $K(X)(q)\to K(X)[q,q^{-1}]$.  

\subsection{Genus-0 reconstruction in terms of the $S$-matrix}
Let us choose a formal power series $\tau:=\tau(Q,\mathbf{t})$ with
coefficients in $K(X)$, such that 
\ben
\left. [S(\tau,q)\mathbf{q}(q)]_+\right|_{q=1} = 0.
\een
Since $[S(\tau,q)(1-q)]_+ = 1-q-\tau$ the above equation is equivalent to
\beq\label{fixed-pt}
\left. [S(\tau,q)\mathbf{t}(q)]_+\right|_{q=1} = \tau.
\eeq
This is a fixed point problem for $\tau$ and we can construct a formal solution using the
standard iteration procedure, e.g.,  put $\tau^{(0)}(Q,\mathbf{t})=0$ and set 
\ben
\tau^{(n+1)}:=
\left. [S(\tau^{(n)},q)\mathbf{t}(q)]_+\right|_{q=1}.
\een  
It is easy to check that as $n\to \infty$, the sequence of formal series
$\tau^{(n)}(Q,\mathbf{t})$ has a limit which provides a solution to
our fixed-point problem.

\begin{lemma}\label{jet-id}
Let $\tau=\tau(Q,\mathbf{t})$ be a solution to the fixed-point
equation \eqref{fixed-pt}. Then  the following formulas hold
\ben
\mathcal{F}(\mathbf{t}) & = &  \frac{1}{2}\langle\mathbf{q}(L),\mathbf{q}(L)
\rangle_{0,2}(\tau),\\
\partial_{n,i}\mathcal{F}(\mathbf{t}) & = &  \langle
\Phi_i(L-1)^n, \mathbf{q}(L)\rangle_{0,2}(\tau),\\
\partial_{n,i}\partial_{m,j}\mathcal{F}(\mathbf{t}) & = &  \langle
\Phi_i(L-1)^n, \Phi_j(L-1)^m \rangle_{0,2}(\tau).
\een
\end{lemma}
\proof
The first formula follows immediately from \eqref{F-anc}, because if
$\tau$ is a solution to \eqref{fixed-pt}, then 
$\overline{\mathcal{F}}_\tau([S(\tau,q)\mathbf{q}(q)]_+)=0$. Indeed,
put $\overline{q}=[S(\tau,q)\mathbf{q}(q)]_+$ and note that the
ancestor correlator
\ben
\langle \overline{\mathbf{q}}(\overline{L}_1),\dots, 
\overline{\mathbf{q}}(\overline{L}_{n})
\rangle_{0,n}(\tau) = \sum_{\ell=0}^\infty \sum_d \frac{Q^d}{\ell!} 
\langle \overline{\mathbf{q}}(\overline{L}_1),\dots, 
\overline{\mathbf{q}}(\overline{L}_{n}),
\tau,\dots,\tau
\rangle_{0,n+\ell}
\een
contains a factor of $(\overline{L}_1-1)\cdots
(\overline{L}_n-1) $, which is a pullback of $(L_1-1)\cdots (L_n-1)=0\in
K(\overline{\mathcal{M}}_{0,n})$.

To prove the second identity, let us differentiate formula \eqref{F-anc} in
$q_{n,i}$ before we specialize $\tau$ to
$\tau(Q,\mathbf{t})$. We get 
\ben
\partial_{n,i} \mathcal{F} &&= \langle \Phi_i
(L-1)^n,\mathbf{q}(L)\rangle_{0,2}(\tau)+ \\
&&+
\sum_{n=0}^\infty\frac{1}{n!} 
\langle
[S(\tau,q)\Phi_i(L-1)^n]_+,
\overline{\mathbf{q}}(\overline{L}_2),\dots, 
\overline{\mathbf{q}}(\overline{L}_{n+1})
\rangle_{0,n+1}(\tau).
\een 
Let us specialize $\tau$ to $\tau(Q,\mathbf{t})$, then
$\overline{\mathbf{q}}(1)=0$ and since $(L_2-1)\cdots (L_{n+1}-1)=0$
in $K(\overline{\mathcal{M}}_{0,n+1})$ we get that the sum on the second
line of the above formula vanishes. 
Note that the same argument can be used to prove the third formula. 
\qed

\subsection{Proof of Theorem \ref{t2}}. We claim that
any solution $\tau=\tau(Q,\mathbf{t})$ to the fixed point problem
\eqref{fixed-pt}  satisfies \eqref{is-alpha} and 
\ben
J(\tau,0)=1+\sum_{j=0}^N\partial_{0,j}\partial_{0,0}\mathcal{F}(\mathbf{t})\,
\Phi^j. 
\een
The above equation follows immediately from Lemma \ref{jet-id} and the
string equation
\ben
&&
\langle 1,\mathbf{t}(L),\dots,\mathbf{t}(L)\rangle_{0,n+1,d} =\\ 
&&
\langle \mathbf{t}(L),\dots,\mathbf{t}(L)\rangle_{0,n,d}+
\sum_{i=1}^n 
\left\langle
\mathbf{t}(L),\dots,
\Delta\mathbf{t}(L),\dots
\mathbf{t}(L)
\right\rangle_{0,n,d}
\een
where the insertion $\Delta\mathbf{t}(L):=
\frac{\mathbf{t}(L)-\mathbf{t}(1)}{L-1}$ is at the i-th place and the
equality holds for all $d\neq 0$ or $n\geq 3$ (see \cite{Lee1,TT}). 
We just have to prove that $\tau=\tau(Q,\mathbf{t})$
is a solution to \eqref{is-alpha}.

We have
\beq\label{J-der}
\partial_{n,i}J(\tau,0) = \partial_{0,0}\, 
\sum_{j=0}^N
\langle \Phi_j,\Phi_i(L-1)^n\rangle_{0,2}(\tau)\,
\Phi^j .
\eeq
Note that if $f\in\CC[q,q^{-1}]$, then 
\beq\label{proj-K_+}
\operatorname{Res}_{q=0,\infty}\, \frac{f(q)dq}{L-q}
= f(L).
\eeq
Therefore
\ben
\langle \Phi_j,\Phi_i(L-1)^n\rangle_{0,2}(\tau)=
\operatorname{Res}_{q=0,\infty}\, 
dq\, (q-1)^n
\left\langle \Phi_j, \frac{\Phi_i}{L-q}\right\rangle_{0,2}(\tau).
\een
By the definition of the operator series $S(\tau,q)$, we get
\ben
\partial_{0,0}\, \left\langle \Phi_j,
  \frac{\Phi_i}{L-q}\right\rangle_{0,2}(\tau)= 
-q^{-1} \partial_{0,0}\,  G(\Phi_j,S(\tau,q)\Phi_i).
\een
Recall that 
\ben
G(a,b) = g(a,S(\tau,0)^{-1}b), 
\een
hence
\ben
\partial_{0,0}\, \left\langle \Phi_j,
  \frac{\Phi_i}{L-q}\right\rangle_{0,2}(\tau)= 
-q^{-1} \partial_{0,0}\,  g(\Phi_j, S(\tau,0)^{-1} S(\tau,q)\Phi_i)
\een
and the equation \eqref{J-der} takes the form
\beq\label{J-der2}
\partial_{n,i}J(\tau,0) = 
-\partial_{0,0}\, 
\operatorname{Res}_{q=0,\infty}\, 
\frac{dq}{q}\, (q-1)^n S(\tau,0)^{-1} S(\tau,q)\Phi_i.
\eeq
Using that $S(\tau,q)$ is a fundamental solution for $\nabla$, we get
\ben
\partial_{n,i}J(\tau,0)=\partial_{n,i}S(\tau,0)^{-1}1 =
S(\tau,0)^{-1} \partial_{n,i}\tau
\een
and 
\ben
\partial_{0,0} (S(\tau,0)^{-1} S(\tau,q))= 
q S(\tau,0)^{-1} \partial_{0,0}S(\tau,q).
\een
The equation \eqref{J-der2} takes the form
\ben
\partial_{n,i}\tau = 
- 
\operatorname{Res}_{q=0,\infty}\, 
dq\, (q-1)^n \partial_{0,0}(S(\tau,q)\Phi_i).
\een
Since $S(\tau,q)$ does not have a pole at $q=0$, the above formula
reduces to the formula that we have to prove.
\qed

\subsection{The Topological recursion relations}

We would like to derive the K-theoretic Topological Recursion Relations (TRR) from the
fact that the solution $\tau=\tau(Q,\mathbf{t})$ to the fixed point
problem \eqref{fixed-pt} is a solution to the system
\eqref{is-alpha}. Put
\ben
\langle 
E_1L_1^{k_1},\dots,E_r L_r^{k_r}\rangle_{0,r}(\mathbf{t}):= 
\sum_{\ell=0}^\infty \sum_d
\frac{Q^d}{\ell!}
\langle \Phi_{i_1}
E_1L_1^{k_1},\dots,E_r L_r^{k_r},
\mathbf{t},\dots,\mathbf{t}
\rangle_{0,r+\ell,d},
\een
where $\mathbf{t}\in K(X)[q,q^{-1}]$ and the $\mathbf{t}$-insertions
in the correlator should be understood as $\mathbf{t}(L_{r+i})$,
$1\leq i\leq \ell$. 
The TRR can be stated as follows.
\begin{proposition}
If $k>0$ and $E_2,E_3\in K(X)$, then
the following identity holds
\ben
&&
\left\langle \Phi_i (L_1-1)^k, E_2L_2^{k_2},E_3L_3^{k_3}\right\rangle_{0,3}(\mathbf{t})
= \\
&&
\sum_{a,b=0}^N 
\left\langle 
\Phi_i L_1(L_1-1)^{k-1}, \Phi_a
\right\rangle_{0,2}(\mathbf{t}) \, 
G^{ab}(\mathbf{t})
\left\langle 
\Phi_b, E_2L_2^{k_2},E_3L_3^{k_3}
\right\rangle_{0,3}(\mathbf{t}),
\een
where $G^{ab}(\mathbf{t})$ are the entries of the matrix inverse to
the matrix $G(\mathbf{t})=(G_{ab}(\mathbf{t}))_{a,b=0}^N $ defined by 
\ben
G_{ab}(\mathbf{t}) := g_{ab}+\langle\Phi_a,\Phi_b\rangle_{0,2}(\mathbf{t}).
\een
\end{proposition}
\proof
According to Lemma \ref{jet-id} we have
\ben
\left\langle E_2L_1^{k_2},E_3L_2^{k_3}\right\rangle_{0,2}(\mathbf{t})
= 
\left\langle E_2L_1^{k_2},E_3L_2^{k_3}\right\rangle_{0,2}(\tau).
\een
Let us write $\mathbf{t}(q) = \sum
t_{k,i}\Phi_i(L-1)^k$. Differentiating the above identity by $t_{k,i}$
we get 
\beq\label{3-jet}
\left\langle \Phi_i (L_1-1)^k,
  E_2L_2^{k_2},E_3L_3^{k_3}\right\rangle_{0,3}(\mathbf{t}) = 
\left\langle
\partial_{k,i}\tau, E_2L_2^{k_2},E_3L_3^{k_3}
\right\rangle_{0,3}(\tau).
\eeq
Recalling that $\tau$ is a solution to 
\eqref{is-alpha} we get 
\ben
\partial_{k,i}\tau=-\operatorname{Res}_{q=0,\infty}
dq (q-1)^{k-1} \partial_{0,0}\tau\bullet S(\tau,q)\Phi_i.
\een
For $k=0$, since $S(\tau,q)$ does not have a pole at $q=0$ or
$\infty$,  the above equation gives
\ben
\partial_{0,i}\tau = \partial_{0,0}\tau\bullet \Phi_i. 
\een
Therefore, we have 
\ben
\partial_{k,i}\tau = 
-\operatorname{Res}_{q=0,\infty}
dq (q-1)^{k-1} 
\sum_{a,b=0}^N G(S(\tau,q)\Phi_i,\Phi_a) G^{ab}(\tau) \partial_{0,b}\tau.
\een
Recalling the definition of the $S$-matrix we get 
\ben
-\operatorname{Res}_{q=0,\infty}
dq (q-1)^{k-1} 
G(S(\tau,q)\Phi_i,\Phi_a) = 
\operatorname{Res}_{q=0,\infty}
dq q(q-1)^{k-1} 
\left\langle \Phi_a,\frac{\Phi_i}{L-q}
\right\rangle_{0,2}(\tau).
\een
The above residue is precisely $\langle \Phi_a,\Phi_i
L(L-1)^{k-1}\rangle_{0,2}(\tau)$ (see formula \eqref{proj-K_+}). The RHS of \eqref{3-jet} turns into 
\ben
\sum_{a,b=0}^N 
\left\langle\Phi_i L_1(L_1-1)^{k-1},\Phi_a\right\rangle_{0,2}(\tau)
G^{ab}(\tau)
\left\langle
\partial_{0,b}\tau, E_2L_2^{k_2},E_3L_3^{k_3}
\right\rangle_{0,3}(\tau).
\een
Note that 
\ben
\left\langle
\partial_{0,b}\tau, E_2L_2^{k_2},E_3L_3^{k_3}
\right\rangle_{0,3}(\tau) = \partial_{0,b} 
\left\langle
E_2L_1^{k_2},E_3L_2^{k_3}
\right\rangle_{0,2}(\tau),
\een
so to complete the proof it remains only to recall the 2-jet identity
in Lemma \ref{jet-id}. 
\qed

\subsection{Geometric proof of TRR }
The TRR can be proved geometrically (see \cite{TT}, Appendix A) using
the standard cotangent lines comparison techniques. Namely, let
$\pi:X_{0,\ell+3,d}\to \overline{\mathcal{M}}_{0,3}$ be the map that
forgets the map to the target $X$ and the last $\ell$ marked
points. Then we have 
\ben
L_1=\overline{L}_1\otimes \mathcal{O}(D),\quad \overline{L}_1=\pi^*L_1,
\een
where $D\subset X_{0,n+3,d}$ is the divisor consisting of stable maps
for which the irreducible component that caries the first marked point
is contructed by $\pi$. Since $\mathcal{O}(-D)=1-\mathcal{O}_D$ and
$\overline{L}_1=1$ we get 
\ben
L_1-1=L_1\otimes \mathcal{O}_D.
\een
The divisor $D$ is a union of divisors $D_i$ obtained from an
appropriate gluing map
\ben
X_{0,\ell'_i+2,d'_i}\times_{\Delta}X_{0, \ell''_i+3,d_i''} \to X_{0,\ell+3,d},
\een
where $\Delta\subset X\times X$ is the diagonal and the fiber product
on the LHS is defined via the evaluation maps 
corresponding to the marked points to be glued. 
The structure sheaf $\mathcal{O}_D $ can be
expressed in terms of the structure sheaves $\mathcal{O}_{D_i}$ using
the inclusion-exclusion principle (see \cite{Gi1})
\ben
\mathcal{O}_D=\sum_{i} \mathcal{O}_{D_i} -\sum_{i,j}
\mathcal{O}_{D_i\cap D_j} + \sum_{i,j,k}
\mathcal{O}_{D_i\cap D_j\cap D_k}-  \cdots. 
\een
Writing the first insertion of the correlator on the LHS of TRR as $\Phi_i (L_1-1)^{k-1}
L_1\otimes \mathcal{O}_D$ and using the above formula for
$\mathcal{O}_D$ we get precisely the RHS of TRR, where
$G^{ab}(\mathbf{t})$ comes as the $(a,b)$-entry of 
\ben
g^{-1} + \sum_{\ell=1}^\infty (-1)^\ell g^{-1}
(G(\mathbf{t})g^{-1})^\ell.
\qed
\een

 Using TRR and Lemma \ref{jet-id} we can give an alternative proof of Theorem \ref{t2}. The
idea is to use TRR to simplify the RHS of formula \eqref{J-der}, which
would allow us to express the derivatives
$\partial_{n,i}\tau(Q,\mathbf{t})$ in terms of the $S$-matrix. The
argument is straightforward, so we skip the details.

\bibliographystyle{amsalpha}

\end{document}